\documentclass{amsart}
\usepackage{amssymb}
\usepackage{amsmath}
\usepackage{graphicx}
\usepackage{amsthm}
\usepackage{dsfont,textcomp}  
\usepackage{mathrsfs}

\begin{document}
\date{April 28th 2006}
\thanks{The author is funded by a Steno Research Grant  from The Danish Natural Science Research Council}
\newtheorem{theorem}{Theorem }[section]
\newtheorem{lemma}[theorem]{Lemma}
\newtheorem{corollary}[theorem]{Corollary}
\newtheorem{proposition}[theorem]{Proposition}
\theoremstyle{definition}
\newtheorem{definition}[theorem]{Definition}
\newtheorem{example}[theorem]{Example}
\theoremstyle{remark}
\newtheorem{remark}[theorem]{Remark}

\renewcommand{\labelenumi}{(\roman{eumi})} 
\def\theenumi{\roman{enumi}}

\numberwithin{equation}{section}

\def \g {{\gamma}}
\def \G {{\Gamma}}\def \l {{\lambda}}
\def \a {{\alpha}}
\def \b {{\beta}}
\def \f {{\phi}}
\def \w {{\omega}}
\def \r {{\rho}}
\def \R {{\mathbb R}}
\def \H {{\mathbb H}}
\def \N {{\mathbb N}}
\def \C {{\mathbb C}}
\def \Z {{\mathbb Z}}
\def \F {{\Phi}}
\def \Q {{\mathbb Q}}
\def \e {{\epsilon }}
\def \GmodH {{\Gamma\backslash\H}}
\def \psl  {{\hbox{PSL}_2( {\mathbb R})} }
\def \vol {{\hbox{vol}}}

\newcommand{\norm}[1]{\left\lVert #1 \right\rVert}
\newcommand{\normm}[1]{\left\lVert #1 \right\rVert_{m}}
\newcommand{\abs}[1]{\left\lvert #1 \right\rvert}

\newcommand{\modsym}[2]{\left \langle #1,#2 \right\rangle}
\newcommand{\inprod}[2]{\left \langle #1,#2 \right\rangle}
\newcommand{\Nz}[1]{\left\lVert #1 \right\rVert_z}
\newcommand{\ver}[1]{\operatorname{vert}\left( #1 \right)}
\newcommand{\wl}[1]{\operatorname{wl}\left( #1 \right)}
\newcommand{\wlc}[1]{\operatorname{wl_c}\left( #1 \right)}
\renewcommand{\labelenumi}{(\roman{enumi})} 
\def\theenumi{\roman{enumi}} 

\title{On pairs of prime geodesics with fixed homology difference. }

\author{Morten S. Risager}
\address{Department of Mathematics, University of Aarhus, Ny Munkegade Building 530, 8000 Aarhus C, Denmark}
\email{risager@imf.au.dk}

\begin{abstract}We exhibit the analogy between prime geodesics on hyperbolic Riemann surfaces and ordinary primes. We present new asymptotic counting results concerning pairs of prime geodesics with fixed  homology difference.
\end{abstract}
\maketitle
\section{Introduction}  
Let $M$ be a compact Riemann surface $M$ of genus $g>1$. It is a
fascinating fact that the norms of the prime closed geodesics on $M$
in many respects are analogous to ordinary primes $p\in \N$. One may
think of them as being \lq pseudo-primes\rq. A striking instance of this analogy is the prime geodesic theorem which was proved by Huber \cite{huber} and  Selberg (see \cite{hejhal}):
\begin{equation}\label{huber-selberg}
\pi(x)=\#\{\g\in\mathscr{P}(M)\vert N(\gamma) \leq x\}\sim li(x).
\end{equation} 
Here $\mathscr{P}(M)$ is the set of prime closed geodesics (a geodesic is prime if it is not an iterate of another geodesic), $li(x)=\int_1^x 1/\log(t)dt$, and $N(\gamma)$ is the norm of $\gamma$ defined by $N(\gamma)=e^{l(\gamma)}$ where $l(\gamma)$ is the geodesic length of $\gamma$.

Consider now $\F:\mathscr{P}(M)\to H_1(M,\Z)$, i.e. the projection to the first homology group with integer coefficients. Fix $\b\in H_1(M,\Z)$ and let $\pi_\b(x)$ be the number of prime geodesics $\gamma$ of norm at most $x$ and with $\F(\gamma)=\b$. Phillips and Sarnak \cite{phillips-sarnak} (and immediately following them Adachi and Sunada \cite{adachi-sunada}) found an asymptotic expansion for $\pi_\beta(x)$:

\begin{equation}\label{phillips-sarnak}
\pi_\b(x)\sim (g-1)^g\frac{x}{\log^{g+1}x}\left(1+\frac{c_1(\beta)}{\log x}+\frac{c_2(\beta)}{\log^2 x}+\cdots \right).
\end{equation} 
The way in which $c_i(\b)$ depends on the specific homology class $\b$ remained unexamined in \cite{phillips-sarnak}. 
We notice that the main term does not depend on $\beta$. 

In certain applications we would like to understand the dependence of
the homology class in this asymptotic expansions. One result in this
direction is the following due to Sharp \cite{sharp2}: Fix an isomorphism $\psi: H_1(M,\Z)\cong \Z^{2g}$
There exist a $2g\times2g$ positive definite symmetric matrix $N$ of determinant 1 such that
\begin{equation}\label{sharp}
\pi_\b(x)=\frac{e^{-\inprod{\psi(\beta)}{N^{-1}\psi(\beta)}/2\sigma^2\log x}}{(2\pi\sigma^2\log x)^g}li(x)+o\left(\frac{x}{\log^{g+1}(x)}\right),
\end{equation}
where $\sigma^{-2}=2\pi(g-1)$ and the implied constant is
\emph{independent} of $\b$. The main point is of course the
independence of $\b$ in the error term, since without this
(\ref{sharp}) reduces to a statement about the main term in a asymptotic expansion a la  (\ref{phillips-sarnak}).

On average we can get better error terms.
Petridis and Risager \cite{petridis-risager} proved that
\begin{equation}\label{petridis-risager}
\sum_{\substack{\b\in B\\\normm{\psi(\b)}\leq \sqrt{\log x}\log\log x}}\left(\pi_\b(x)-\frac{e^{-\inprod{\psi(\beta)}{N^{-1}\psi(\beta)}/2\sigma^2\log x}}{(2\pi\sigma^2\log x)^g}li(x)\right)=o(li(x)).
\end{equation}

From
(\ref{petridis-risager}) follows an equidistribution result concerning geodesics in (large) sets of homology classes: 
 Let $\norm{r}$ be
a norm on $\R^{2g}$. We say that a subset $B\subseteq H_1(M,\Z)$ has asymptotic density $d_{\norm{\cdot}}(B)$ with respect to $\norm{\cdot}$ if the limit
\begin{equation*}\lim_{x\to\infty}\frac{\#\{\b\in B\vert \norm{\psi(\b)}\leq
  x\}}{\#\{\b\in H_1(M,\Z)\vert \norm{\psi(\b)}\leq x\}}\end{equation*}
exists and equals $d_{\norm{\cdot}}(B)$, i.e. if the image in $\Z^{2g}$ under $\psi$
has asymptotic density in $\Z^{2g}$  with respect to $\norm{\cdot}$.
In \cite{petridis-risager} it was shown that there exist a norm $\norm{\cdot}_M$ such that for all sets $B\subseteq H_1(M,\Z)$ with asymptotic density with respect to $\norm{\cdot}_M$
\begin{equation}
\frac{\pi_B(x)}{\pi(x)}\to d_{\norm{\cdot}_M}(B)\textrm{ as } x\to\infty.
\end{equation}
Here $\pi_B(x)$ is the number of prime closed geodesics with norm $N(\gamma)$ at
most $x$ and homology class $\F(\g)\in B$.

One may investigate what happens if we consider pairs  -- or
more generally $k$-tuples -- of prime closed geodesics. 
Pollicott and Sharp \cite{sharp-pollicot} recently did so in the following way: Let $a_1,\ldots,a_g,b_1,\ldots b_g$ be a fundamental
set of generators for the fundamental group $\pi_1(M)$ (see section
\ref{generalsection}) The conjugacy classes of $\pi_1(M)$ are in one to one
correspondence with closed geodesics on $M$. For a closed geodesic we
let $\abs{\g}=\min\{\wl{g} g\in\{\g \}\}$ where
$\{\g\}$ is the conjugacy class associated with the closed
geodesic $\g$ and $\wl\g$ is the word length of $g$ in the fundamental
set of generators. Pollicott and Sharp used sub-shifts of finite
type and the thermodynamic formalism to prove the following pair
correlation result:  there exist a
constant $c$ such that for any $a<b$
\begin{equation}
  \label{sharp-pollicot}
\#\{(\g,\g')\vert \abs{\g}, \abs{\g'}\leq n, \, a\leq
l(\g)-l(\g')\leq b\}\sim c(b-a)\frac{e^{2n}}{n^{5/2}},   
\end{equation}
in the limit $n\to \infty$. We notice that in terms of the \lq
pseudo-primes\rq{}  $N(\gamma)$ Pollicott and Sharp are looking at
\emph{quotients of norms} in finite intervals. For a somewhat different type
of results  concerning pairs see \cite{pitt}.

In this paper we study a new - more geometric - counting functions for pairs of geodesics. More precisely we consider the
counting function for pairs of prime closed geodesics with norm at most
$x$ and fixed homology difference:
\begin{equation}
   \label{pair-def}
    \pi_2^\beta(x)=\#\{\g_1,\g_2\in \mathscr{P}(M)\vert N(\g_i)\leq x,
    \, \F(\g_2)-\F(\g_1)=\beta\}.
  \end{equation}
This counting function is geometric in the sense that the ordering of elements is according to the geodesic length.
We will prove the following result:
\enlargethispage{\baselineskip}
\begin{theorem}\label{pairs-asymp}Let $\b\in H_1(M,\Z)$.
  \begin{equation*}
    \pi_2^\b(x)\sim\frac{(g-1)^g}{2^g}\frac{x^2}{\log^{g+2}(x)}
  \end{equation*}
in the limit $x\to\infty$.
\end{theorem}
In particular there are infinitely many pairs of prime geodesics with fixed homology difference.
One may think of this as a hyperbolic Riemann surface version of the twin prime conjecture. For
further explanation as to this analogy we refer to section
\ref{ordinaryprimes}. We can now ask how the error term depends on the
specific homology class $\beta$. We prove the following result:

\begin{theorem}\label{mylocallimit} Let $\b\in H_1(M,\Z)$.
  \begin{equation*}
    \pi_2^\b(x)=\frac{1}{2^g}\frac{e^{-\inprod{\psi(\beta)}{N^{-1}\psi(\beta)}/2\sigma^2\log(x)}}{(2\pi\sigma^2\log(x))^g}\frac{x^2}{\log^{2}(x)}+o\left(\frac{x^2}{\log^{g+2}(x)}\right)
  \end{equation*}
when $x>3$, where, when $\normm{\psi(\b)}=o(\sqrt{\log x}/\log \log x)$,  the implied constant is independent of $\beta$.
\end{theorem}
As with Sharps result (\ref{sharp}) the main point in Theorem
\ref{mylocallimit} is the existence of an error term which is
independent of $\beta$.  Theorem \ref{pairs-asymp} follows trivially
from Theorem \ref{mylocallimit}. The assumption  $\normm{\psi(\b)}=o(\sqrt{\log x}/\log \log x)$ may be relaxed to $\normm{\psi(\b)}=o(\sqrt{\log x})$. 

\begin{remark}
The geometry of the surface $M$ is intimately linked with the spectrum
of the Laplacian of the surface considered as a Riemannian
manifold. This link is evident from the Selberg trace formulae which
relates the lengths of closed geodesics with the eigenvalues of the
Laplacian in a summation formulae. (See (\ref{stf}) below). This
\lq duality\rq{}  between the length spectrum and the Laplace spectrum
has proven itself extremely useful both in the study of
eigenvalues (e.g. Weyl's law (see e.g. \cite[\S
4.4]{venkov})) as well as in the study of the lengths of geodesics
which is what we investigate in the present work. We use the Selberg trace formulae to count primes in a  homology class (a technique developed by Phillips and Sarnak \cite{phillips-sarnak}), and we keep track of the dependence on the specific homology class in the error terms. We then analyze how these error terms contribute to the relevant sum. 
\end{remark}

\begin{remark}The fact that we are considering surfaces of fixed
  negative sectional curvature $-1$, is not essential. If $M$ has
  variable negative curvature we can combine the ideas of this paper
  with the ideas developed by Sharp \cite{sharp2}, to get results
  similar to theorems \ref{pairs-asymp} and \ref{mylocallimit}. In
  this case the proof uses the thermodynamic formalism instead of the
  Selberg trace formulae. 
It is also possible to obtain similar results for free groups using
ideas by Petridis and Risager
\cite{petridis-risager,petridis-risager-free-groups}.    

\end{remark}

\begin{remark}
The techniques used in this paper may be used to study counting results of a more general type than the ones we consider. Consider any $A\subseteq H_1(M,\Z)^k$. We may then consider the counting function
\begin{equation}\label{moregeneral}
  \#\{(\g_i)\in  \mathscr{P}(M)^k\vert N(\g_i)\leq x,
    \,(\phi(\g_i))\in A \}
\end{equation}
This may be rewritten as 
\begin{equation*}
  \sum_{(\a_i)_{i=1}^k\in A}\prod_{i=1}^k\pi_{\a_i}(x)
\end{equation*}
By using good expansions for  $\pi_{\a_i}(x)$ it is now possible to study (\ref{moregeneral}). To prove Theorem \ref{mylocallimit} we develop and analyze this in full for 
\begin{equation*}
  A=\{(\a_1,\a_2)\in H_1(M,\Z)^2\vert \a_2-\a_1=\beta\}
\end{equation*}
 The techniques certainly apply to much more general sets. 
We hope this paper will be facilitating for anyone interested in such questions.
\end{remark}
The paper is organized as follows: In section 2 we describe how the results in this introduction may be seen as analogues of statements in analytic number theory. Section 3 briefly describe the technique developed by Phillips and Sarnak (combined by an idea of Sharp \cite{sharp2}) to count primes in a specific homology class. In the following section we describe how this may be transformed into counting prime pairs with fixed homology difference. In Section 5 we find the main and error term of a counting function with certain weights, and in Section 6 we explain how to use multi-variable summation by parts to remove these weights. 

\section{Ordinary primes in arithmetic progressions}\label{ordinaryprimes}
There is nothing new in this section. Its purpose is  to
emphasize how (almost) all the results mentioned in the introduction
are analogues of classical results or conjectures in analytic number
theory. Readers not interested in such connections should feel free to
move to the next section, as the rest of the paper does not depend
directly on this section. We quote from \cite{iwaniec-kowalski}
but most of the results can be found in any solid textbooks on analytic number theory. 

Let $\Pi(x)=\#\{p\leq x\}$ be the number of primes less than or equal to
$x$. The prime number theorem \cite[Section 2.1]{iwaniec-kowalski} proved by Hadamard and de la Vall\'ee Poussin asserts that 
\begin{equation}\label{pnt}
\Pi(x)\sim li(x).
\end{equation}
The theorem of Huber and Selberg (\ref{huber-selberg}) is analogous to (\ref{pnt}).  

Given a primitive conjugacy class $a \bmod q$ i.e. $(a,q)=1$ we let
$\Pi(x;a,q)$ be the number of primes less than $x$ with $p\equiv a \bmod q$. The main result about primes in arithmetic progressions is (\cite[(17.2)]{iwaniec-kowalski})

\begin{equation}\label{progressions}
\Pi(x;q,a)\sim \frac{li(x)}{\Phi(q)},
\end{equation}
where $\Phi(q)=\#\{1\leq a< q\vert (a,q)=1\}$ is the Euler totient. The result of Phillips and Sarnak (\ref{phillips-sarnak}) my be considered analogous to (\ref{progressions}).

In applications to other problems involving primes it is  of great
interest to know how the error term in (\ref{progressions}) depends on
$q$ and $x$. A first result in this direction is the Siegel-Walfisz
theorem \cite[Corollary 5.29]{iwaniec-kowalski} which states that for
all $A>0$ $a, q\in\N$, $(a,q)=1$  
\begin{equation}
\Pi(x;q,a)= \frac{li(x)}{\Psi(q)}+ O\left(\frac{x}{log^A(x)}\right),
\end{equation}
when $x>2$ where the implied constant \emph{depends only} on $A$.
 We like to think of Sharps theorem (\ref{sharp}) as analogous to this result.

The extremely  potent idea of taking averages over $q$ and $a$ to get better bounds on the error term on average has been used very successfully in the famous theorem of Bombieri and Vinogradov \cite[Theorem 17.1]{iwaniec-kowalski}: 

\begin{theorem}[Bombieri-Vinogradov] \label{bomb-vino} For any $A>0$ there exist $B>0$ such that
\begin{equation*}
\sum_{q\leq Q}\max_{(a,q)=1}\abs{\Pi(x;q,a)-\frac{li(x)}{\Phi(q)}}=O\left(\frac{x}{\log^A(x)}\right) 
\end{equation*}
where $Q=x^{1/2}log^{-B}(x)$. The implied constant depends only on $A$.
\end {theorem}
Conjecturally  (Elliot-Halberstram) we can take $Q=x^{1-\e}$. In many applications Theorem \ref{bomb-vino} is an excellent substitute for the Generalized Riemann hypothesis which says that $\Pi(x,q,a)=li(x)/\Psi(q)+O(x^{1/2+\e})$  The large range $q\leq x^{1-\e}$  can be handled on average if we allow averages in $a$ also (See \cite[Theorem 17.2]{iwaniec-kowalski}):
\begin{theorem}[Barbon, Davenport, Halberstram]\label{barb-dave-halb}
For any $A>0$ there exist $B>0$ such that 
\begin{equation*}
\sum_{q\leq Q}\sum_{\substack{a \bmod q\\(a,q)=1}}\left(\Pi(x;q,a)-\frac{li(x)}{\Phi(q)}\right)^2 =O\left(\frac{x}{\log^A(x)}\right) 
\end{equation*}
where $Q=xlog^{-B}(x)$. The implied constant depends only on $A$.
\end{theorem}
The theorem of Petridis and Risager may be considered analogous to Theorems \ref{bomb-vino} and \ref{barb-dave-halb}.

A folklore conjecture says that there are infinitely many twin primes
i.e. primes $p$ such that $p+2$ is a prime. This conjecture was quantified by Hardy and Littlewood who conjectured that 
\begin{equation}\label{twin-conjecture}
\#\{p_1,p_2\leq x\vert p_2-p_1=2 \}\sim 2c_2\int_1^x\frac{1}{\log^2(t)}dt
\end{equation}
where $c_2=\prod_{p>2}(1-(p-1)^{-2})$. We could prove it if we
where able to handle certain linear combinations of
$\Pi(x;q,a)-li(x)/\Psi(q)$ See (\cite[Section
13.1]{iwaniec-kowalski}). Certainly the Montgomery conjecture --
$\Pi(x;q,a)=li(x)/\Psi(q)+O(x^{1/2+\e}/q^{1/2})$ -- would give it
immediately. Unfortunately we are not able to handle the relevant
linear combinations and the twin prime conjecture remains completely open.

Theorem \ref{pairs-asymp} is analogous to the Conjecture
(\ref{twin-conjecture}), and its proof goes along the same lines as what
one would like to do for primes. But for prime geodesics the
$\beta$ dependence of $\pi_\beta(x)$ can be understood well enough
that we  can prove which contributions
give error terms and which contribution gives the main term in the
relevant linear combination.

\section{Counting prime closed geodesics in homology classes}\label{generalsection}

In this section we set up some notation and explain how the Selberg trace formula can be used to count geodesics in a homology class. We then quote an equality from Petridis and Risager \cite{petridis-risager} which is proved using this technique. This equality is the starting point of our current investigation.

Any compact Riemann surface of genus $g>1$ without
boundary may be realized as $M=\GmodH$, where $\H$ is the upper
half-plane and $\G\subseteq \psl$ is a strictly hyperbolic discrete subgroup of $\psl$ acting on
$\H$ by linear fractional transformations. The surface $M$ has fundamental group $\pi_1(M)=\Gamma$. The
closed oriented geodesics  are in one to one correspondence with the conjugacy
classes of $\G$ by the following recipe: Pick a base point $z_0\in
\H$ above $m$. From a conjugacy class $\{\g \}$ we project
($\bmod\,  \G$)
the geodesic in $\H$ from $z_0$ to $\g z_0$ to $M$ which is homologous to  a
closed geodesic on $M$. 

The group $\G$ has a fundamental set of generators i.e. a set of
generators \begin{equation*}a_1,\ldots,a_g,b_1,\ldots,b_g \in \G\end{equation*} with one defining
relation
\begin{equation*}
  [a_1,b_1]\cdots [a_b,b_g]=1
\end{equation*}
where $[a,b]$ is the commutator of $a$ and $b$.

We let $C_i$, $i=1,\ldots g$, be the geodesic induced by $a_i$  and
$C_{g+i}$, $i=1,\ldots g$, be the geodesic induced by $b_i$ (The
fundamental generators lie in different conjugacy classes so $i\neq j$
implies $C_i\neq C_j$). The first homology group $H_1(M,\Z)$ is
isomorphic to the free group over $\Z$ of $C_1,\ldots,C_{2g}$, i.e.
\begin{equation*}
  H_1(M,\Z)\cong\left\{\sum m_iC_i\vert m_i\in \Z \right\}\cong\Z^{2g}
\end{equation*}

The exist a basis for the space of harmonic
$1$-forms which is dual to $C_1,\ldots C_{2g}$ in the sense that
\begin{equation}
  \int_{C_i}\omega_j=\delta_{ij}.
\end{equation}
These lift to harmonic differentials $\alpha_i=\Re(f_i(z)dz)$ on $\H$ where
$f_i(z)$ is a holomorphic form of weight $2$ with respect to
$\G$. Then $\g\in\G$ induces a geodesic with homology $\sum m_iC_i$ if
and only if
\begin{equation}
  \phi(\g):=\left(\int_{z_0}^{\gamma z_0}\a_1,\ldots, \int_{z_0}^{\gamma z_0}\a_{2g}\right)=(m_1,\ldots,m_{2g}).
\end{equation}
We notice that $\phi(\g)$ does not depend on the choice of path or of the choice of $z_0$. 
Consider the unitary characters on $\G$ defined by
 \begin{equation}\begin{array}{llccc}\chi_\e&:&\G&\to& S^1\\
 &&\g&\mapsto&\displaystyle e^{2\pi i \inprod{\phi(\g)}{\e}}
 \end{array}.
 \end{equation} 
where $\e\in\R^{2g}$ and $\langle\cdot,\cdot\rangle$ is the usual
inner product on $\R^n$.

Consider now the set of the  set of square-integrable $\chi_\e$-automorphic functions, i.e. the set of $f:\H\to\C$ such that
 \begin{equation}\label{automorphic}f(\g z)=\chi_\e(\g)f(z)
 \end{equation}
 and 
 \begin{equation}\label{summable}\int_{F}\abs{f( z)}^2d\mu(z)<\infty ,
 \end{equation} where $F$ is a fundamental domain for $\GmodH$. 
 Let $L_\e$ denote the Laplacian defined as the closure of 
 \begin{equation}
   -y^2\left(\frac{\partial^2}{\partial x^2}+\frac{\partial^2}{\partial y^2}\right)
 \end{equation} defined on smooth compactly supported functions
 satisfying (\ref{automorphic}) and (\ref{summable}). The Laplacian is
 self-adjoint and its spectrum consists of a countable set of
 eigenvalues $0\leq \l_0(\e)\leq \l_1(\e)\leq \ldots$
We write $\lambda_j(\e)=1/4+r_j^2(\e)=s_j(\e)(1-s_j(\e))$. All our
geodesic counting results has their origin in the Selberg trace formula
for $L_\e$ which relates the Laplace spectrum with the length
spectrum in a very precise way. (See \cite{selberg, hejhal}):
 \begin{align}\label{stf}
\begin{aligned}
   \sum_j\hat h(r_j(\e))=&2(g-1)\int_{-\infty}^\infty r\tanh(\pi r)\hat
   h(r)dr\\
&+\sum_{\{\g\}}\frac{\chi_\e(\g)l(\g )}{k\sinh(l(\g)/2)}h(l(\g ))
\end{aligned}
\end{align}
where $h$ is a smooth even function on $\R$ of compact support, $\hat
h$ is its Fourier transform and $l(\g)$ is the length of the geodesic
induced by $\g\in \G$.  When $\e=0$ the contribution from
$r_0(0)$ should be counted twice. One of the main ideas in
\cite{phillips-sarnak} is that by multiplying (\ref{stf}) with
$\exp(-2\pi i \inprod{\psi(\b)}{\e})$ and then integrating over
the whole character variety (i.e. over $\e\in\R^{2g}\backslash
\Z^{2g}$) we pick out exactly those $\gamma$ on the right hand side of
(\ref{stf}) with homology class $\beta$.
 
By combining ideas of Sharp \cite{sharp2} and Phillips and Sarnak
\cite{phillips-sarnak} it is possible 
to get precise information from (\ref{stf}) about 
\begin{equation}
  R_\beta(x)=\sideset{}{'}\sum_{\substack{N(\gamma)\leq
      x\\\F(\g)=\beta}}\frac{l(\g)}{\sinh{(l(\g)/2)}}
\end{equation}
(the ${}^\prime$ on the sum means that we only sum over prime geodesics).
Petridis and Risager noticed \cite[(2.11)]{petridis-risager} that  up to an error term
of  decay (independent of $\beta$)  $x^{-\delta}$ 
\begin{equation}\label{letshope1}
  \frac{R_\b(x)}{4\sqrt{x}}-\frac{e^{-\inprod{\psi(\b)}{N^{-1}\psi(\beta)}/2\sigma^2\log(x)}}{(2\pi \sigma^2\log (x))^g}
\end{equation}
equals 
\begin{equation}\label{letshope2}
  \int_{B(\rho)}
\left(
\frac{e^{(s_0(\e)-1)\log(x)}}{2s_0(\e)-1
}-e^{-\inprod{\e}{N\e}4\pi^2\sigma^2\log(x)/2}
\right)
\overline{\chi^\b_{\e}}d\e
\end{equation}
 for every sufficiently small $\rho$. This will be the
starting point for our investigation concerning pairs of prime geodesic. 
\section{Counting prime pairs with fixed homology difference}
In this section we explain how to use the counting technique described in the previous section to count \emph{pairs} of geodesics with restrictions on their homology difference.

We define, for $x_1,x_2>1$, 
\begin{equation*}
  \pi_2^\b(x_1,x_2):=\#\left\{\g_1,\g_2\in\mathscr{P}(M)\left\vert
      N(\g_i)\leq x_i, \F(\g_2)-\F(\g_2)=\b\right.\right \}
\end{equation*}
and we denote $\pi_2^\b(x):=\pi_2^\b(x,x)$. We fix $0<k<1$. \emph{We will
always assume that}  
\begin{equation}\label{samesize}
  x^k\leq x_i\leq x. 
\end{equation}
An obvious choice is to let $x=\max x_i$. Then (\ref{samesize}) puts
restrictions on $\min x_i$. 
The restriction (\ref{samesize}) implies that $\log(x_1)$, $\log(x_2)$, and $\log(x)$ are all of
the same size (i.e. $\log(x_1)\asymp\log(x_2)\asymp\log(x)$) The same
is true for $\log^{-1}(x_1)$, $\log^{-1}(x_2)$, and
$\log^{-1}(x)$. When we, in the following,  estimate various sums the error term may depend on $k$ but
never on $x$.
 
Instead of working with $\pi_2^\b(x_1,x_2)$ directly it turns out to be more
convenient for us to work with something closer related to
$R_\b(x)$. In principle we would like to use
\begin{align*}
  \pi_2^\b(x_1,x_2)=&\sum_{\a\in
    H_1(M,\Z)}\#\left\{\g_1,\g_2\in\mathscr{P}(M)\left\vert\begin{array}{l} N(\g_i)\leq x_i\\ (\F(\g_1),\F(\g_2))=(\a,\b+\a)\end{array}\right.\right \}\\
=&\sum_{\a\in H_1(M,\Z)}\pi_\a(x_1)\pi_{\b+\a}(x_2)
\end{align*}
but it turns out to be  more convenient to use
\begin{align}
\nonumber  R_2^\beta(x_1,x_2):=&\sideset{}{'}\sum_{\substack{N(\gamma_i)\leq
      x_i\\\F(\g_2)-\F(\g_1)=\beta}}\frac{l(\g_1)l(\g_2)}{\sinh{(l(\g_1)/2)}\sinh{(l(\g_2)/2)}}\\
\label{thisisit}=&\sum_{\a\in H_1(M,\Z)}R_{\a}(x_1)R_{\a+\b}(x_2)
\end{align}

The main strategy is now to use (\ref{letshope1}) and
(\ref{letshope2}) to find an asymptotic expansion for (\ref{thisisit}) and then use multi-dimensional partial summation to get
the expansion for $\pi_2^\b(x_1,x_2)$.

We start by making some  estimates on $R_\beta(x)$.
Consider  
\begin{equation*}f_x(\e)=\left(
\frac{e^{(s_0(\e)-1)\log(x)}}{2s_0(\e)-1
}-e^{-\inprod{\e}{N\e}4\pi^2\sigma^2\log(x)/2}
\right).\end{equation*}
We let 
\begin{eqnarray}\label{itstomuch}
\begin{aligned}
  A(\b,x)&=&4\sqrt{x}\frac{e^{-\inprod{\psi(\b)}{N^{-1}\psi(\b)}/2\sigma^2\log(x)}}{(2\pi \sigma^2\log (x))^g}\\
  B(\b,x)&=&4\sqrt{x}\int_{B(\rho)}
f_x(\e)
\overline{\chi^\b_{\e}}d\e.
\end{aligned}
\end{eqnarray}
From (\ref{letshope1}) and (\ref{letshope2}) we have 
\begin{equation}
  R_\b(x) =A(\b,x)+B(\b,x)+O(x^{1/2-\delta})
\end{equation}
for some $\delta>0$. The constant $\delta$ and the implied constant
are absolute. 

To be able to bound expressions involving $B(\b,x)$, i.e. $\Sigma_2,\Sigma_3,$ and $\Sigma_4$, we recall Proposition 2.5  from \cite{petridis-risager}.

 \begin{proposition}\label{againvanishinlimit}
 Let $N=\{\inprod{\w_i}{\w_j}\}.$ 

 \begin{enumerate}
 \item{\label{againfirst}For every
 $\e_0\in \R^{2g}$ 
  \begin{equation*}
     e^{(s_0(\e/2\pi\sigma\sqrt{\log(x)})-1)\log(x)}\to
       e^{-\inprod{\e}{N\e}/2}
   \end{equation*}
 as $x\to\infty$.}
 \item{\label{againsecond}There exists $\delta>0$ such that for all $\norm{\e}<\delta\sqrt{\log(x)}$ 
   \begin{equation*}
     \abs{e^{(s_0(\e/2\pi \sigma\sqrt{\log(x)})-1)\log(x)}-
       e^{-\inprod{\e}{N\e}/2}}\leq 2e^{-\inprod{\e}{N\e}/4}.
   \end{equation*}}
 \item{\label{againthirde}For all $\theta>0$ sufficiently small there exist  $ C>0$
     such that for all  $\log(x)>0$, $\norm{\e}<\delta \log(x)^{\theta}$, 
 \begin{equation*}
     \abs{e^{(s_0(\e/\r\sqrt{\log(x)})-1)\log(x)}-
       e^{-\inprod{\e}{N\e}/2}}\leq C\frac{1}{\log(x)^{1-2\theta}}.
 \end{equation*}}
  \item{\label{againforth} Let $0< \nu<1/4$. For every $k>0$ there exist constants  $\delta_1,\delta_2>0$
     such that, 
 \begin{equation*}
     \abs{e^{(s_0(\e/2\pi \sigma\sqrt{\log(x)})-1)\log(x)}-
       e^{-\inprod{\e}{N\e}/2}}\leq \frac{e^{-\nu\inprod{\e}{N\e}}}{\log^k(x)}.
 \end{equation*}}
 when $\log(x)>0$, $\delta_1\sqrt{\log(\log (x))}<\norm{\e}<\delta_2 \sqrt{\log(x)}$.
\end{enumerate}
 \end{proposition}

We use Proposition \ref{againvanishinlimit} to prove the following lemma:
\begin{lemma}\label{detherblirgrimt}Let $j=1,2$.  For every $l>0$ there exist a $\delta>0$ such that 
  \begin{align}\label{endnuendnuen} \int_{B'(x)} \abs{f_x(\e)}^j d\e&=O((\log x)^{-(g+j)+\varepsilon})\\ \label{noken}
\int_{B(\delta)\backslash B'(x)}\abs{f_x(\e)}^jd\e&=O((\log x)^{-l})\end{align}
where  $B'(x)=B(\delta\sqrt{\log\log x}/\sqrt{\log x})$.   

\end{lemma}
\begin{proof}
We let $j=1$. By a change of variables we see that 
\begin{align*}
\int_{B'(x)}\abs{f_x(\e)}d\e&=(2\pi\sigma\sqrt{\log x})^{-2g}\!\!\!\!\!\!\!\!\!\!\!\!\!\int\limits_{\norm{\e}\leq 2\pi \sigma\delta \sqrt{\log \log x}}\abs{f_x(\e/2\pi \sigma\sqrt{\log x})}d\e\\
  \int_{B(\delta)\backslash B'(x)}\abs{f_x(\e)}d\e&=(2\pi\sigma\sqrt{\log x})^{-2g}\!\!\!\!\!\!\!\!\!\!\!\!\!\!\!\!\!\!\!\!\!\!\int\limits_{2\pi \sigma\delta \sqrt{\log \log x}\leq \norm {\e}\leq \rho 2\pi \sigma\sqrt{\log (x)}}\!\!\!\!\!\!\!\!\!\!\!\!\!\!\!\!\!\!\abs{f_x(\e/2\pi \sigma\sqrt{\log x})}d\e
\end{align*}

We now let $\delta=\delta_1/2\pi\sigma$ where $\delta_i$ are constants as in Proposition \ref{againvanishinlimit} (\ref{againforth}) with $k=l$.
 We may safely assume that $\rho$ has been chosen so small that $\rho2\pi\sigma<\delta_2$.
 Since $s_0(\e)$  is even with $s_0(0)=1$ we have 
 \begin{equation}\label{dinnerwasheavy}
   \abs{(2s_0(\e)-1)^{-1}-1}\leq C \norm{\e}^2
 \end{equation}
when $\norm{\e}\leq \rho$. Hence the integrand is 
bounded  by 
\begin{equation*}
  \abs{e^{(s_0(\e/2\pi\sigma\sqrt{\log x})-1)\log x}-e^{-\inprod{\e}{N\e}/2}}+C\abs{e^{(s_0(\e/2\pi\sigma\sqrt{\log x})-1)\log x}}\norm{\e}^2\log x^{-1},  
\end{equation*}
which by Proposition \ref{againvanishinlimit} (\ref{againsecond}) is bounded by \begin{equation*}
  \abs{e^{(s_0(\e/2\pi\sigma\sqrt{T})-1)T}-e^{-\inprod{\e}{N\e}/2}}+Ce^{-\mu\inprod{\e}{N\e}}\log x^{-1},   
\end{equation*}
for some small $\mu>0$. 

When $\norm{\e}\leq 2\pi \sigma\delta \sqrt{\log \log x}$ we use Proposition \ref{againvanishinlimit} (\ref{againthirde}) to conclude (\ref{endnuendnuen}).

We can safely assume that $\delta_2$ is big enough that 
\begin{equation*}
e^{-\mu\inprod{\e}{N\e}}\leq \frac{e^{-\mu/2\inprod{\e}{N\e}}}{\log ^l(x)}
\end{equation*} when $\delta_2 \sqrt{\log \log x}\leq \norm {\e}\leq \rho 2\pi \sigma\sqrt{\log (x)}$. Using Proposition \ref{againvanishinlimit} (\ref{againthirde}) we see that 
\begin{equation*}
  \abs{f_x(\e/2\pi \sigma\sqrt{\log x})}\leq C' \frac{e^{-\mu\inprod{\e}{N\e}/2}}{\log ^l(x)}
\end{equation*} in this region, from which we easily conclude (\ref{noken}). 

The case $j=2$ is similar.

\end{proof}

We now return to the sum (\ref{thisisit}) We start by
showing that we only need a finite sum:
\begin{lemma}\label{boundedlemma}Let $\normm{r}=\sum\abs{r_i}$ be the max norm. There exist a constant  $C>0$ depending only on $M$  such
  that 
  \begin{equation}
    \normm{\phi(\g)}\leq C l_\g
  \end{equation}
for all closed geodesics $\g$ .
\end{lemma}
\begin{proof}
  This follows directly from \cite[Lemma 2.4]{petridis-risager}.
\end{proof}

From Lemma \ref{boundedlemma} follows that $R_\b(x)=0$ if $\normm{\psi(\b)}> C
\log(x)$. This implies that in (\ref{thisisit}) we only need to sum over 
\begin{align}
  \label{trivial1}
  \normm{\psi(\a+\b)}&\leq C\log x_2\\ \label{trivial2} \normm{\psi(\a)}&\leq C\log x_1
\end{align}
Since we are mainly interested in asymptotics we may restrict the sum to a much smaller sum. We define $u(x)=\sqrt{\log(x)}\log\log x$ and let 
\begin{equation}
  \tilde R_2^\beta(x_1,x_2):=\sum_{\substack{\a\in H_1(M,\Z)\\\normm{\psi(\a)}\leq u(x)}}R_{\a}(x_1)R_{\a+\b}(x_2)
\end{equation}
Then 
\begin{align}
  R_2^\beta(x_1,x_2)-\tilde R_2^\beta(x_1,x_2)&=\sum_{\substack{\a\in H_1(M,\Z)\\u(x)\leq \normm{\psi(\a)}\leq C\log (x)}}R_{\a}(x_1)R_{\a+\b}(x_2)\\
&\leq O\Biggl(\frac{\sqrt{x_2}}{\log^g x_2}\sum_{\substack{\a\in
    H_1(M,\Z)\\u(x)\leq \normm{\psi(\a)}\leq C\log
    (x)}}R_{\a}(x_1)\Biggr)  
\end{align}
We used Lemma \ref{detherblirgrimt} to get the inequality.  The sum
\begin{equation*}
  \sum_{\substack{\a\in H_1(M,\Z)\\u(x)\leq \normm{\psi(\a)}\leq C\log (x)}}R_{\a}(x_1)
\end{equation*}
is $o(\sqrt{x_1})$ by \cite[Lemma 2.7]{petridis-risager}.

It follows that 
\begin{equation}\label{smallerisenough}
  R_2^\beta(x_1,x_2)-\tilde R_2^\beta(x_1,x_2)=o\left(\frac{\sqrt{x_1x_2}}{\log^g(x)}\right)
\end{equation}
so for the asymptotic results we are aiming at,  we may consider the small sum $\tilde R_2^\beta(x_1,x_2)$.

\enlargethispage{\baselineskip}
Using (\ref{smallerisenough}) we conclude that when $x\geq 3$
\begin{align}
\nonumber  R_2^\b(x_1,x_2)={}&\sum_{\substack{\a\in H_1(M,\Z)\\ \normm{\psi(\a)}\leq
u(x)}}A(\a,x_1)A(\a+\b,x_2)\\
\nonumber&+\sum_{\substack{\a\in H_1(M,\Z)\\ \normm{\psi(\a)}\leq
u(x)}}A(\a,x_1)B(\a+\b,x_2)
\label{uglyone}\\&+\sum_{\substack{\a\in H_1(M,\Z)\\ \normm{\psi(\a)}\leq
u(x)}}B(\a,x_1)A(\a+\b,x_2)\\
\nonumber&+\sum_{\substack{\a\in H_1(M,\Z)\\ \normm{\psi(\a)}\leq
u(x)}}B(\a,x_1)B(\a+\b,x_2) + o((x_1x_2)^{1/2}/\log^g x)\\
\nonumber={}&\Sigma_1(\b,x_1,x_2)+\Sigma_2(\b,x_1,x_2)+\Sigma_3(\b,x_1,x_2)+\Sigma_4(\b,x_1,x_2)\\
&\qquad+o((x_1x_2)^{1/2}/\log^g x).
\nonumber 
\end{align}
for some $\delta'>0$.
 
\section{Finding the main term and error term}
We ended the last section by splitting the function $R_2^\b(x_1,x_2)$ into four different contributions and an error term. In this section we determine which contributions are \lq big\rq{} and which are \lq small\rq{}. We will prove the following result:
\begin{theorem}\label{main-R} Let $\b\in H_1(M,\Z)$ and $0<k<1$. Then 
  \begin{equation*}
    R_2^{\b}(x_1,x_2)=16x_1^{1/2}x_2^{1/2}\frac{e^{-\inprod{\psi(\beta)}{N^{-1}\psi(\beta)}/2\sigma^2\log x_2}}{(2\pi\sigma^2(\log x_1+\log x_2))^g}+o\left(\frac{x_1^{1/2}x_2^{1/2}}{\log^{g/2}(x_1)\log^{g/2}(x_2)}\right)
  \end{equation*}
when $3<x^k\leq x_i\leq x$, 
$\normm{\psi(\b)}=o(\sqrt{\log x}/\log\log x))$ and $x\to\infty$,  where the implied
constant depends at most on $k$ and $M$.
\end{theorem}
Our starting point is the identity (\ref{uglyone}). We start by showing that the main term comes out of $\Sigma_1$.
\begin{lemma}\label{ug1}
  \begin{equation*}
    \Sigma_1(\beta,x_1,x_2)=16x_1^{1/2}x_2^{1/2}\frac{e^{-\inprod{\psi(\beta)}{N^{-1}\psi(\beta)}/2\sigma^2\log x_2}}{(2\pi\sigma^2(\log x_1+\log x_2))^g}+o\left(\frac{x_1^{1/2}x_2^{1/2}}{\log^{g/2}(x_1)\log^{g/2}(x_2)}\right)
  \end{equation*}
where, when $\normm{\psi(\b)}=o(\sqrt{\log x}/\log\log x)$, the implied constant is independent of~$\beta$.
\end{lemma}
\begin{proof}
  Using (\ref{uglyone}) and (\ref{itstomuch}) we easily find that
  \begin{align*}
\Sigma_1(\beta,x_1,x_2)=16\frac{\sqrt{x_1x_2}}{(2\pi\sigma^2)^{2g}}&\frac{1}{(\log
  x_1\log x_2)^g} 
\\&\cdot\sum_{\substack{\a\in \Z^{2g}\\\normm{\a}\leq u(x)}}e^{-\frac{\inprod{\a}{N^{-1}\a}}{2\sigma^2\log x_1}} e^{-\frac{\inprod{(\a+\psi(\b))}{N^{-1}(\a+\psi(\b))}}{2\sigma^2\log x_2}}\\
=16\frac{\sqrt{x_1x_2}}{(2\pi\sigma^2)^{2g}}&\frac{e^{-\inprod{\psi(\beta)}{N^{-1}\psi(\beta)}/2\sigma^2\log
    x_2}}{(\log x_1+\log x_2)^g}
\\&\cdot \sum_{\substack{\a\in \Z^{2g}\\\normm{\a}\leq u(x)}}\frac{e^{-\frac{\inprod{\a}{N^{-1}\a}}{2\sigma^2g(x_1,x_2)}}}{g(x_1,x_2)^g} e^{-2\frac{\inprod{\a}{N^{-1}\psi(\b)}}{2\sigma^2\log x_2}}
 \end{align*}

where \begin{equation}g(x_1,x_2)=(\log^{-1}x_1+\log^{-1}x_2)^{-1}=\frac{\log x_1\log x_2}{\log x_1+\log x_2}.\end{equation} We must therefore understand the above sum.

We consider first the case where $\b=0$. In this case we consider
\begin{equation}
  \frac{1}{(2\pi\sigma^2)^g}\sum_{\substack{\a\in \Z^{2g}\\\normm{\a}\leq u(x)}}\frac{e^{-\inprod{\a}{N^{-1}\a}/2\sigma^2g(x_1,x_2)}}{g(x_1,x_2)^g}
\end{equation}
It follows from \cite[Lemma 2.10]{petridis-risager} that this sum converges to 1, and we conclude that
\begin{equation*}\Sigma_1(0,x_1,x_2)=\frac{16x_1^{1/2}x_2^{1/2}}{(2\pi\sigma^2(\log x_1+\log x_2))^g}+o\left(\frac{x_1^{1/2}x_2^{1/2}}{\log^{g/2}(x_1)\log^{g/2}(x_2)}\right)\end{equation*} which proves the lemma when $\beta=0$. 

The general case follows in the same way  if we verify that 
\begin{align}
\label{lastmanstanding}  \sum_{\substack{\a\in \Z^{2g}\\\normm{\a}\leq u(x)}}\frac{e^{-\inprod{\a}{N^{-1}\a}/2\sigma^2g(x_1,x_2)}}{g(x_1,x_2)^g}&(1-e^{-2\inprod{\a}{N^{-1}\psi(\b)}/2\sigma^2\log x_2})=o\left(1\right)
\end{align}
when $\normm{\beta}=o(\sqrt{\log x}/\log\log x)$.

There exist a decreasing function $r(x)$ going to zero as $x\to\infty$ such that $\normm{\beta}\leq r(x)\sqrt{\log x}/\log \log x).$ Hence there exist an absolute constant $C>0$

\begin{equation}(1-e^{-2\inprod{\a}{N^{-1}\psi(\b)}/2\sigma^2\log x_2})\leq C r(x)\end{equation} 
when $\normm{\a}\leq \sqrt{\log x}\log \log x$. The bound (\ref{lastmanstanding}) now follows from this and the fact that the remaining part of the sum converges to 1 (\cite[Lemma 2.10]{petridis-risager}). In particular it is bounded (independently of $\beta$) \end{proof}

\begin{remark}
It is fairly straightforward to improve Lemma \ref{ug1} to only requiring $\normm{\beta}=o(\sqrt{\log x}).$ This requires showing that we need only to sum over $u(t)=\sqrt{\log (x)}v(x)$  in \ref{uglyone}. Here   $v(x)$ is some function which grows sufficiently slowly to infinity. 
\end{remark}

We are now ready to bound the 3 remaining terms of (\ref{uglyone}) which all go into the error term in Theorem \ref{main-R}:

\begin{lemma}\label{ug2}
  \begin{equation*}
    \Sigma_i(\beta,x_1,x_2)=o\left(\frac{x_1^{1/2}x_2^{1/2}}{\log^{g/2}(x_1)\log^{g/2}(x_2)}\right),\qquad i=2,3
  \end{equation*}
where the implied constant is independent of $\beta$.
\end{lemma}
\begin{proof}
We see from (\ref{itstomuch}) and (\ref{uglyone}) that the function $\Sigma_3(\b,x_1,x_2)$ equals
\begin{equation*}
 16x_1^{1/2}x_2^{1/2}\int\limits_{B(\rho)}f_{x_1}(\e)\!\!\!\!\!\!\!\!\!\sum_{\substack{\a\in H_1(M,\Z)\\ \normm{\psi(\a)}\leq
u(x)}}\!\!\!\!\!\!\!\!\!\!\!\frac{e^{-\inprod{\psi(\a+\b)}{N^{-1}\psi(\a+\b)}/2\sigma^2\log(x_2)}}{(2\pi
  \sigma^2\log (x_2))^g} \overline{\chi^\a_{\e}}d\e.
\end{equation*}
By bounding the character trivially and using Lemma \ref{detherblirgrimt} we get that
\begin{equation}\label{hereren}
  \Sigma_3(\beta,x_1,x_2)=O\Biggl(\frac{x_1^{1/2}x_2^{1/2}}{\log^{g+1-
      \varepsilon}(x_1)}  
  \sum_{\substack{\a\in H_1(M,\Z)\\ \normm{\psi(\a)}\leq 
u(x)}} 
\mkern-18mu
\frac{e^{-\inprod{\psi(\a+\b)}{N^{-1}\psi(\a+\b)}/2\sigma^2\log(x_2)}}{(2\pi 
  \sigma^2\log (x_2))^g}   \Biggr).  
\end{equation}
The result follows from the fact that 
\begin{equation*}\log^{-1}(x_1)\leq k^{-1}\log^{-1}(x_2),\end{equation*} if we show that the remaining sum in   (\ref{hereren})  is uniformly bounded.
To see this we note that 
\begin{align*}
  \sum_{\substack{\a\in H_1(M,\Z)\\ \normm{\psi(\a)}\leq
u(x)}}\!\!\!\!\!\!\!\!\!\!\!&\frac{e^{-\inprod{\psi(\a+\b)}{N^{-1}\psi(\a+\beta)}/2\sigma^2\log(x_2)}}{(2\pi
  \sigma^2\log (x_2))^g}  \leq  \sum_{\substack{\a\in H_1(M,\Z)\\ \normm{\psi(\a)}\leq
u(x)+\normm{\beta} }}\!\!\!\!\!\!\!\!\!\!\!\frac{e^{-\inprod{\psi(\a)}{N^{-1}\psi(\a)}/2\sigma^2\log(x_2)}}{(2\pi
  \sigma^2\log (x_2))^g} \\
&\leq \sum_{\substack{\a\in H_1(M,\Z)\\ \normm{\psi(\a)}\leq
u(x) }}\!\!\!\!\!\!\!\!\!\!\!\frac{e^{-\inprod{\psi(\a)}{N^{-1}\psi(\a)}/2\sigma^2\log(x_2)}}{(2\pi
  \sigma^2\log (x_2))^g}+\!\!\!\!\!\!\!\!\!\!\! \sum_{\substack{\a\in H_1(M,\Z)\\ u(x)\leq\normm{\psi(\a)}\leq
u(x)+\normm{\beta} }}\!\!\!\!\!\!\!\!\!\!\!\frac{e^{-\inprod{\psi(\a)}{N^{-1}\psi(\a)}/2\sigma^2\log(x_2)}}{(2\pi
  \sigma^2\log (x_2))^g} \\
\end{align*}
 The first sum is independent of $\beta$ and converging to 1 (\cite[Lemma 2.10]{petridis-risager}. The second sum may be bounded by an absolute constant times
 \begin{align*}
  \sum_{\substack{n \in \N\\ u(x)\leq n}} \frac{e^{-\mu n^2/\log x}}{
  \sqrt{\log x}}  \left( \sum_{n \in \N} \frac{e^{-\mu n^2/\log x}}{
  \sqrt{\log x}} \right)^{g-1}
 \end{align*}
for some small $\mu>0.$ It is easy to see that this is bounded (in fact converging to zero) by comparison with an integral.
 The sum $\Sigma_2(\beta,x)$ may be  handled in the
same way.
\end{proof}

\begin{lemma}\label{ug3}
  \begin{equation*}
    \Sigma_4(\beta,x_1,x_2)=o\left(\frac{x_1^{1/2}x_2^{1/2}}{\log^{g/2}(x_1)\log^{g/2}(x_2)}\right)
  \end{equation*}
where the implied constant is independent of $\beta$.
\end{lemma}
\begin{proof}
  Using the definitions (\ref{itstomuch}) and (\ref{uglyone}) we see
  that $\Sigma_4(\b,x_1,x_2)$ equals
  \begin{equation}\label{youknowwhotalkstomuch}
    16x_1^{1/2}x_2^{1/2}\!\!\!\int\limits_{B(\rho)\times B(\rho)}f_{x_1}(\e^{1})f_{x_2}(\e^{2})\sum_{\substack{\a\in H_1(M,\Z)\\ \normm{\psi(\a)}\leq
u(x)}} \overline{\chi^\a_{\e^1}\chi^{\a+\b}_{\e^2}}d\e^{1}
d\e^{2}
  \end{equation}
We write 
$B(\rho)\times B(\rho)=[B'(x)\times B'(x)]\cup [(B(\rho)\times B(\rho))\setminus (B'(x)\times B'(x))]$ where $B'(x)$ is defined in Lemma \ref{detherblirgrimt}, and the corresponding $l$ is chosen sufficiently large. We split the integral in (\ref{youknowwhotalkstomuch}) accordingly. 

In the part over $B(\rho)\times B(\rho)\backslash B'(x)\times B'(x)$, we bound the sum trivially by $O(u(v)^{2g})$ and the remaining integral is easily seen to be $O(\log(x)^{-h})$ for $h$ as big as we care for by choosing the data of $B'(x)$ sufficiently well.

In the integral over $B'(x)\times B'(x)$ we use Cauchy-Schwarz to conclude that 
\begin{equation*}    \int\limits_{B'(x)\times B'(x)}f_{x_1}(\e^{1})f_{x_2}(\e^{2})\sum_{\substack{\a\in H_1(M,\Z)\\ \normm{\psi(\a)}\leq
u(x)}} \overline{\chi^\a_{\e^1}\chi^{\a+\b}_{\e^2}}d\e^{1}
d\e^{2}
  \end{equation*}
is bounded by 
\begin{equation}
\Biggl(\:\int\limits_{B'(x)\times
  B'(x)}\abs{f_{x_1}(\e^{1})f_{x_2}(\e^{2})}^2d\e^{1} 
d\e^{2}\Biggr)^{1/2},
\end{equation}
 which is $O((\log x_1\log x_2)^{-(g/2+1)+\varepsilon})$ by Lemma
 \ref{detherblirgrimt},  times  
\begin{equation}
  \Biggl(\:\int\limits_{B'(x)\times B'(x)}\abs{\sum_{\substack{\a\in H_1(M,\Z)\\ \normm{\psi(\a)}\leq
u(x)}} \overline{\chi^\a_{\e^1}\chi^{\a+\b}_{\e^2}}}^2d\e^{1}
d\e^{2}\Biggr)^{1/2}
\end{equation}

Notice that 
\begin{align*}
  \Biggl(\:\int\limits_{B'(x)\times B'(x)} &  \abs{\sum_{\substack{\a\in H_1(M,\Z) \\ \normm{\psi(\a)}\leq
u(x)}} \overline{\chi^\a_{\e^1}\chi^{\a+\b}_{\e^2}}}^2d\e^{1}
d\e^{2}\Biggr)^{1/2}
\\
&\leq \#\{\a\in H_1(M,\Z) |\, \normm{\psi(\a)}\leq
u(x)\}\times \vol(B'(x))\\ &=O(u(x)^{2g}(\sqrt{\log\log x}/\sqrt{\log x})^{2g})\\
&=O((\log\log x)^{3g})
\end{align*}
Collecting the pieces we see that 
 \begin{equation*}
    \Sigma_4(\beta,x_1,x_2)=O\left(\frac{x_1^{1/2}x_2^{1/2}}{\log^{g/2}(x_1)\log^{g/2}(x_2)}\log^{\varepsilon-1} x\right)
  \end{equation*}
 from which the result follows.
\end{proof}
Theorem \ref{main-R} now follows from  (\ref{uglyone}), Lemmata  \ref{ug1},
\ref{ug2}, and \ref{ug3}.

\section{Using partial summation}
In this section we show how to use multi-dimensional partial summation to conclude from
Theorem \ref{main-R} our
main result: 
\begin{theorem}\label{main} Let $\b\in H_1(M,\Z)$ and $0< k<1$.
  \begin{align*}
    \pi_2^\b(x_1,x_2)=\frac{e^{-\inprod{\psi(\beta)}{N^{-1}\psi(\beta)}/2\sigma^2\log(x_2)}}{(2\pi\sigma^2(\log x_1+\log(x_2))^g}&\frac{x_1x_2}{\log(x_1)\log(x_2)}\\ &+o\left(\frac{x_1x_2}{\log^{g/2+1}(x_1)\log^{g/2+1}(x_2)}\right) 
  \end{align*}
when $3<x^k\leq x_i\leq x$ and $x\to\infty$, where the implied
constant depends at most on $k$ and $M$ for $\normm{\psi(\b)}=o(\sqrt{\log x}/\log \log x)$.
\end{theorem}
We notice that putting $x_1=x_2$ we obtain Theorem \ref{mylocallimit} and Theorem \ref{pairs-asymp}

To prove Theorem \ref{main} we let \begin{equation*}P_2^\b(x_1,x_2)=\sideset{}{'}\sum_{\substack{N(\gamma_i)\leq
      x_i\\\F(\g_2)-\F(\g_1)=\beta}}\frac{4\log N(\g_1)\log
    N(\g_2)}{\sqrt{N(\g_1)}\sqrt{N(\g_2)}}).\end{equation*}
We have a trivial bound
\begin{equation}\label{thisistrivial}
P_2^{\b}(x_1,x_2)=O(x_1^{1/2}x_2^{1/2})  
\end{equation}
which follows directly from (\ref{huber-selberg}) by ignoring the
condition $\F(\g_2)-\F(\g_1)=\beta$.

It is not difficult to see that
\begin{align}\nonumber&\abs{R_2^\beta(x_1,x_2)-P_2^\beta(x_1,x_2)}\\
\nonumber&\qquad\leq \sideset{}{'}\sum_{N(\gamma_i)\leq x_i}
\frac{2\log(N(\g_1))}{\sqrt{N(\gamma_1)}}\frac{\log(N(\g_2))}{N(\g_2)\sinh{(\log(N(\g_2)/2))}}\\
\nonumber&\qquad\quad-\sideset{}{'}\sum_{N(\gamma_i)\leq x_i}
\frac{2\log(N(\g_2))}{\sinh(\log(N(\g_2)/2))}
\frac{\log(N(\g_1))}{N(\g_1)\sinh{(\log(N(\g_1)/2))}}\\
&\qquad=\sideset{}{'}\sum_{N(\gamma_1)\leq x_1}\label{easier} 
\frac{2\log(N(\g_1))}{\sqrt{N(\gamma_1)}}\sideset{}{'}\sum_{N(\gamma_2)\leq x_2}\frac{\log(N(\g_2))}{N(\g_2)\sinh{(\log(N(\g_2)/2))}}\\
\nonumber&\qquad\quad-\sideset{}{'}\sum_{N(\gamma_2)\leq x_2}
\frac{2\log(N(\g_2))}{\sinh(\log(N(\g_2)/2))}
\sideset{}{'}\sum_{N(\gamma_1)\leq x_1}\frac{\log(N(\g_1))}{N(\g_1)\sinh{(\log(N(\g_1)/2))}}\\
\nonumber&\qquad=O(x_1^{1/2}+x_2^{1/2})
\end{align}
We used (\ref{huber-selberg}) again in the last estimate. It follows that
Theorem \ref{main-R} holds with $R_2^{\b}(x_1,x_2)$ replaced by $P_2^{\b}(x_1,x_2)$.

Using multi-dimensional partial summation (\cite[Theorem
1.6]{kraetzel})we find 
\begin{align}
 \nonumber \sideset{}{'}\sum_{\substack{N(\gamma_i)\leq
      x_i\\\F(\g_2)-\F(\g_1)=\beta}}1=&\sideset{}{'}\sum_{\substack{N(\gamma_i)\leq
      x_i\\\F(\g_2)-\F(\g_1)=\beta}}\frac{4\log N(\g_1)\log
    N(\g_2)}{\sqrt{N(\g_1)}\sqrt{N(\g_2)}}\cdot\frac{\sqrt{N(\g_1)}\sqrt{N(\g_2)}}{4\log N(\g_1)\log
    N(\g_2)}\\
\nonumber=&\frac{\sqrt{x_1}\sqrt{x_2}}{4\log x_1\log x_2}P_2^\b(x_1,x_2)\\
\nonumber&\quad-\frac{\sqrt{x_1}}{4\log
  x_1}\int_1^{x_2}P_2^\b(x_1,t_2)m(t_2)dt_2\\
\label{gettingcloser}&\quad-\frac{\sqrt{x_2}}{4\log
  x_2}\int_1^{x_1}P_2^\b(t_1,x_2)m(t_1)dt_1\\
\nonumber&\quad+\frac{1}{4}\int_1^{x_1}\int_1^{x_2}P_2^\b(t_1,t_2)m(t_1)m(t_2)dt_2dt_1.
\end{align}
where \begin{equation}m(t)=\frac{d}{dt}\frac{\sqrt{t}}{\log(t)}=\frac{1}{2\sqrt{t}\log(t)}-\frac{1}{\sqrt{t}\log^2(t)}\end{equation}
We then find the asymptotics of the three integrals.

\begin{lemma}\label{stupidlemma1}Let $w\in\R$.
\begin{align*}
  \int_1^{x_1}P_2^\b(t_1,x_2)
  (\sqrt{t_1}\log^w(t_1))^{-1}dt_1=\frac{\sqrt{x_1}}{\log^w(x_1)}&P_2^\b(x_1,x_2)\\ &+o\left(\frac{x_1x_2^{1/2}}{\log^{g/2+w}(x_1)\log^{g/2}(x_2)}\right)
\end{align*}
where the implied constant is independent of $\b$.
\end{lemma}
\begin{proof}
We start by noting that it follows from the trivial bound (\ref{thisistrivial}) that
 \begin{equation*}\int_1^{\sqrt{x_1}}P_2^\b(t_1,x_2)
  (\sqrt{t_1}\log^w(t_1))^{-1}dt_1=o(x_1^{3/4}x_2^{1/2})\end{equation*}
(Clearly  $P_2^\b(x_1,x_2)=0$ for $x_1$ or $x_2$ close to 1 since
lengths of geodesics does not accumulate at zero, so the integral
makes sense even though there seems to be a singularity at 1). 

Let \begin{equation*}m_\b(x_1,x_2)=16x_1^{1/2}x_2^{1/2}\frac{e^{-\inprod{\psi(\beta)}{N^{-1}\psi(\beta)}/2\sigma^2\log x_2}}{(2\pi\sigma^2(\log x_1+\log x_2))^g}\end{equation*}
Consider now 
\begin{align}
\nonumber  \int_{\sqrt{x_1}}^{x_1}&P_2^\b(t_1,x_2)
  (\sqrt{t_1}\log^w(t_1))^{-1}dt_1-\frac{\sqrt{x_1}}{\log^w(x_1)}P_2^\b(x_1,x_2)\\
\nonumber&=\int_{\sqrt{x_1}}^{x_1}(P_2^\b(t_1,x_2)-m_\b(t_1,x_2))
  (\sqrt{t_1}\log^w(t_1))^{-1}dt_1\\
\label{insertone}&\quad+\int_{\sqrt{x_1}}^{x_1}m_\b(t_1,x_2)(\sqrt{t_1}\log^w(t_1))^{-1}dt_1-\frac{\sqrt{x_1}}{\log^w(x_1)}P_2^\b(x_1,x_2)
\\
\nonumber&=\int_{\sqrt{x_1}}^{x_1}(P_2^\b(t_1,x_2)-m_\b(t_1,x_2))
  (\sqrt{t_1}\log^w(t_1))^{-1}dt_1\\
\nonumber&\quad+\frac{\sqrt{x_1}}{\log^{w}(x_1)}m_\b(x_1,x_2)-\frac{\sqrt{x_1}}{\log^w(x_1)}P_2^\b(x_1,x_2)
+O\left(\frac{x_1^{1/2} x_2^{1/2}}{\log^{g+w+1}(x_1)}\right).
\end{align}

  Using Theorem \ref{main-R} and (\ref{easier}) it follows that there
  exist a function $g_k(x)$ depending on $k'$, independent of $\b$,  and
  decreasing to zero as $x \to
  \infty$ such that if $x^{k'}\leq x_i \leq x$ then 
  \begin{equation}\label{nicetrick}
    \abs{P_2^{\b}(x_1,x_2)-m_\b(x_1,x_2)}\leq g_{k'}(x)\frac{x_1^{1/2}x_2^{1/2}}{\log^{g/2}(x_1)\log^{g/2}(x_2)}   
  \end{equation}

We let $k'=k/2$. Then since we are assuming  $x^{k}\leq x_i \leq x$
we have $x^{k'}\leq t_1,x_2 \leq x$ when $t_1\geq \sqrt{x_1}$.

If we take absolute values in (\ref{insertone}) we therefore find
\begin{align*}
  \abs{\,\cdot\,}&\leq g_{k'}(\sqrt{x})x_2^{1/2}\log^{-g/2}(x_2)\int_{\sqrt{x_1}}^{x_1}\log^{-g/2-w}(t_1))dt_1\\
&\quad+\frac{\sqrt{x_1}}{\log^{w}(x_1)}g_{k'}(x)\frac{x_1^{1/2}x_2^{1/2}}{\log^{g/2}(x_1)\log^{g/2}(x_2)}   
+O\left(\frac{x_1 x_2^{1/2}}{\log^{g+w+1}(x_1)}\right)\\
&=o\left(\frac{x_1x_2^{1/2}}{\log^{g/2+w}(x_1)\log^{g/2}(x_2)} \right).
\end{align*}
\end{proof}
\begin{lemma}\label{stupidlemma2}Let $w\in\R$.
\begin{align*}
  \int_1^{x_2}P_2^\b(x_1,t_2)(\sqrt{t_2}\log^w(t_2))^{-1}dt_2=\frac{\sqrt{x_2}}{\log^w(x_2)}&P_2^\b(x_1,x_2)\\&+o\left(\frac{x_1^{1/2}x_2}{\log^{g/2}(x_1)\log^{g/2+w}(x_2)}\right)
\end{align*}
where the implied constant is independent of $\b$.
\end{lemma}
\begin{proof}
We claim that 
\begin{equation*}
  \int_{\sqrt{x_2}}^{x_2}m_\b(x_1,t_2)(\sqrt{t_2}\log^w(t_2))^{-1}dt=m_\b(x_1,x_2)\frac{\sqrt{x_2}}{\log^w(x_2)}+O\left(\frac{x_1^{1/2}x_2}{\log^{g/2}(x_1)\log^{g/2+w+1}}\right)
\end{equation*}
where the implied constant depends at most on $k$ and $w$.  
  Using partial integration we see that
  \begin{align}
\nonumber \int_{\sqrt{x_2}}^{x_2}&\frac{e^{-\inprod{\psi(\b)}{N^{-1}\psi(\b)}/2\sigma^2\log(t_2)}}{\log^{g+w}(t_2)}dt_2=x_2\frac{e^{-\inprod{\psi(\b)}{N^{-1}\psi(\b)}/2\sigma^2\log(x_2)}}{\log^{g+w}(x_2)}\\
&\nonumber\quad-\sqrt{x_2}\frac{e^{-\inprod{\psi(\b)}{N^{-1}\psi(\b)}/2\sigma^2\log(\sqrt{x_2})}}{\log^{g+w}(\sqrt{x_2})}  \\  
\label{intm}&\quad-\int_{\sqrt{x_2}}^{x_2}\frac{e^{-\inprod{\psi(\b)}{N^{-1}\psi(\b)}/2\sigma^2\log(t_2)}(-\inprod{\psi(\b)}{N^{-1}\psi(\b)}/2\sigma^2\log(t_2)-(g+w))}{\log^{g+w+1}(t_2)}dt_2.
  \end{align}
Since $e^{-av}v$ is bounded for $v\in\R_+$ the enumerator of the
integrand is bounded (depending on $g$ an $w$) and the claim follows
easily.

Using the claim the proof follows the proof of Lemma \ref{stupidlemma1} almost verbatim.
\end{proof}
\begin{lemma}\label{stupidlemma3}Let $w_1,w_2\in\R$.
  \begin{align*}\int_1^{x_1}\!\!\!\int_1^{x_2}P_2^\b(t_1,t_2)&(\sqrt{t_1}\log^{w_1}(t_1)\sqrt{t_2}\log^{w_2}(t_2))^{-1}dt_1dt_2\\
&=\frac{x_1^{1/2}x_2^{1/2}}{\log^{w_1}(x_1)\log^{w_2}(x_2)}P_2^\b(x_1,x_2)\\
&\quad +o\left(\frac{x_1x_2}{\log^{g/2+w_1}(x_1)\log^{g/2+w_2}(x_2)}\right)
\end{align*}
\end{lemma}
\begin{proof}
  If we bound $P_2^\b(x_1,x_2)$  trivially \ref{thisistrivial} we see that we only
  need to bound the integral over $(t_1,t_2)\in [\sqrt{x_1},x_1]\times
  [\sqrt{x_2},x_2]$.
\begin{align*}&\abs{\int_{\sqrt{x_1}}^{x_1}\int_{\sqrt{x_2}}^{x_2}\cdots
    dt_1dt_2-\frac{x_1^{1/2}x_2^{1/2}}{\log^{w_1}(x_1)\log^{w_2}(x_2)}P_2^\b(x_1,x_2)} 
\\
& \quad
\leq\abs{\int_{\sqrt{x_1}}^{x_1}\int_{\sqrt{x_2}}^{x_2}(P_2^\b(t_1,t_2)-m_\b(t_1,t_2))(\sqrt{t_1}\log^{w_1}(t_1)\sqrt{t_2}\log^{w_2}(t_2))^{-1}dt_1dt_2}
\\
&\qquad+\abs{\int_{\sqrt{x_1}}^{x_1}\int_{\sqrt{x_2}}^{x_2}
  \frac{m_\b(t_1,t_2)}{\sqrt{t_1}
    \log^{w_1}(t_1)\sqrt{t_2}\log^{w_2}(t_2)}dt_1dt_2-
  \frac{x_1^{1/2}x_2^{1/2}}{\log^{w_1}(x_1) \log^{w_2}(x_2)}P_2^\b(x_1,x_2)} 
\\
\intertext{We use (\ref{intm}) and a calculation on the last integral:}
& \quad
=\abs{\int_{\sqrt{x_1}}^{x_1}\int_{\sqrt{x_2}}^{x_2}(P_2^\b(t_1,t_2)-m_\b(t_1,t_2))(\sqrt{t_1}\log^{w_1}(t_1)
  \sqrt{t_2}\log^{w_2}(t_2))^{-1}dt_1dt_2}
\\ 
&\qquad+\abs{\frac{x_1x_2}{\log^{w_1}(x_1)\log^{w_2}(x_2)}(m_\b(x_1,x_2)-P_2^\b(x_1,x_2))}+
O\left(\frac{x_1x_1}{\log^{g+w_1+w_2}(x)}\right)
\\ 
\intertext{We then use (\ref{nicetrick})}
&\quad=g_{k'}(x)\int_{\sqrt{x_1}}^{x_1}\int_{\sqrt{x_2}}^{x_2}\frac{t_1^{1/2}t_2^{1/2}}{\log^{g/2}(t_1)
  \log^{g/2}(t_2)}(\sqrt{t_1}\log^{w_1}(t_1)\sqrt{t_2}\log^{w_2}(t_2))^{-1}dt_1dt_2
\\ 
&\qquad+g_{k'}(x)\frac{x_1^{1/2}x_2^{1/2}}{\log^{w_1}(x_1)\log^{w_2}(x_2)}
\frac{x_1^{1/2}x_2^{1/2}}{\log^{g/2}(x_1)\log^{g/2}(x_2)}
+O\left(\frac{x_1x_1}{\log^{g+w_1+w_2}(x)}\right)\\ 
&\quad =o\left(\frac{x_1x_2}{\log^{g/2+w_1}(x_1)\log^{g/2+w_2}(x_2)}\right).
\end{align*}
which finishes the proof the the lemma.
\end{proof}
We are now ready to finish the proof of Theorem \ref{main}. From
(\ref{gettingcloser}) and lemmata
\ref{stupidlemma1}, \ref{stupidlemma2} and \ref{stupidlemma3} we
find that 
\begin{align}
 \nonumber \pi_2^\b(x_1,x_2)=\frac{1}{16}\frac{x_1^{1/2}x_2^{1/2}}{\log x_1\log
    x_2}&P_2^\b(x_1,x_2)\\
&+O\left(\frac{x_1^{1/2}x_2^{1/2}}{\log x_1\log
    x_2\log(x)}P_2^\b(x_1,x_2)\right)\\&+o\left(\frac{x_1x_2}{\log^{g/2+1}(x_1)\log^{g/2+1}(x_2)}\right)\nonumber
\end{align}
From (\ref{gettingcloser}) and Theorem \ref{main-R} we find easily
\begin{equation}
P_2^\b(x_1,x_2)=0\left(\frac{x_1^{1/2}x_2^{1/2}}{\log^g(x)}\right)
\end{equation}
where the implied constant is independent of $\b$. We conclude that 
\begin{equation}\label{almostthere}
\pi_2^\b(x_1,x_2)=\frac{1}{16}\frac{x_1^{1/2}x_2^{1/2}}{\log(x_1)\log(x_2)}P_2^\b(x_1,x_2)+o\left(\frac{x_1x_2}{\log^{g/2+1}(x_1)\
\log^{g/2+1}(x_2)}\right)
\end{equation}
Using (\ref{almostthere}) Theorem \ref{main} follows from Theorem \ref{main-R}.

\end{document}